\documentclass[pdflatex,sn-nature]{sn-jnl}

\usepackage{graphicx}
\usepackage{multirow}
\usepackage{amsmath,amssymb,amsfonts}
\usepackage{amsthm}
\usepackage[title]{appendix}
\usepackage{xcolor}
\usepackage{textcomp}
\usepackage{manyfoot}
\usepackage{booktabs}
\usepackage{listings}
\usepackage{array}
\usepackage{tabularx}
\usepackage{placeins}
\usepackage{siunitx}
\newcommand{\code}[1]{\texttt{#1}}
\theoremstyle{thmstyleone}

\theoremstyle{thmstyletwo}

\theoremstyle{thmstylethree}

\usepackage{etoolbox}
\makeatletter
\patchcmd{\@maketitle}{Corresponding author(s). E-mail(s):}{Corresponding author. E-mail:}{}{}
\def\email#1{\global\advance\emailcnt by 1\relax\g@addto@macro\corrauthemail{\textcolor{blue}{#1}}}
\makeatother
\usepackage{pdfpages}
\begin{document}
\newcommand{\ExternalCodexPassCount}{64}
\newcommand{\ExternalClaudePassCount}{64}
\newcommand{\ExternalRecordsPerAgent}{64}
\newcommand{\ExternalCodexPassRate}{100.0\%}
\newcommand{\ExternalClaudePassRate}{100.0\%}
\newcommand{\ExternalCodexRepeatAgreement}{32/32}
\newcommand{\ExternalClaudeRepeatAgreement}{32/32}
\newcommand{\ExternalEndpointReclassifications}{0}
\newcommand{\ExternalPassDelta}{+0.0}
\newcommand{\ExternalPassDeltaLow}{+0.0}
\newcommand{\ExternalPassDeltaHigh}{+0.0}
\newcommand{\ExternalCodexHeterogeneousPassCount}{16}
\newcommand{\ExternalCodexAdvectionPassCount}{16}
\newcommand{\ExternalCodexLShapePassCount}{16}
\newcommand{\ExternalCodexCoupledPassCount}{16}
\newcommand{\ExternalClaudeHeterogeneousPassCount}{16}
\newcommand{\ExternalClaudeAdvectionPassCount}{16}
\newcommand{\ExternalClaudeLShapePassCount}{16}
\newcommand{\ExternalClaudeCoupledPassCount}{16}

\newcommand{\NMIPrimaryScheduledCount}{1450}
\newcommand{\NMIEndpointPassCount}{1162}
\newcommand{\NMIEndpointPassRate}{80.14\%}
\newcommand{\NMIConfirmedRouteValidCount}{1037}
\newcommand{\NMIStrictRouteValidLowerRate}{71.52\%}
\newcommand{\NMIAuditConditionedEstablishedRouteValidCount}{1037}
\newcommand{\NMIAuditConditionedEstablishedRouteValidRate}{71.52\%}
\newcommand{\NMIConfirmedDirectCount}{120}
\newcommand{\NMIConfirmedMixedCount}{5}
\newcommand{\NMIConfirmedRouteInvalidCount}{125}
\newcommand{\NMIConfirmedInvalidRate}{8.62\%}
\newcommand{\NMIRouteUnresolvedCount}{0}
\newcommand{\NMIRouteUnresolvedRate}{0.00\%}
\newcommand{\NMIKnownInvalidCorrectedCount}{1037}
\newcommand{\NMIKnownInvalidCorrectedUpperRate}{71.52\%}
\newcommand{\NMIRouteValidInterval}{71.52--71.52\%}
\newcommand{\NMIAuditConditionedPartialIdentificationInterval}{71.52--71.52\%}
\newcommand{\NMIConfirmedInvalidInflationPoints}{8.62}
\newcommand{\NMIEndpointStrictGapPoints}{8.62}
\newcommand{\NMIConfirmedInvalidEndpointShare}{10.76\%}
\newcommand{\NMIAuditCorpusRecords}{13000}
\newcommand{\NMIAuditCorpusEndpointPasses}{9007}
\newcommand{\NMIAuditCorpusRouteValid}{8882}
\newcommand{\NMIAuditCorpusDirect}{120}
\newcommand{\NMIAuditCorpusMixed}{5}
\newcommand{\NMIAuditCorpusClaimInvalid}{125}
\newcommand{\NMIAuditCorpusUnresolved}{0}
\newcommand{\NMICompletedDoubleReviewCount}{250}
\newcommand{\NMICompletedTriggerReviewCount}{125}
\newcommand{\NMICompletedControlReviewCount}{125}
\newcommand{\NMIFourClassAgreement}{96.0\%}
\newcommand{\NMIFourClassKappa}{0.925}
\newcommand{\NMIBinaryAgreement}{96.8\%}
\newcommand{\NMIBinaryKappa}{0.938}
\newcommand{\NMIAdjudicatedDisagreementCount}{10}
\newcommand{\NMICompletedConsensusNumericalCount}{127}
\newcommand{\NMICompletedConsensusDirectCount}{120}
\newcommand{\NMICompletedConsensusMixedCount}{3}
\newcommand{\NMICompletedConsensusRouteInvalidCount}{123}
\newcommand{\NMIExpandedReviewItemCount}{218}
\newcommand{\NMIExpandedReviewUnresolvedCount}{90}
\newcommand{\NMIExpandedReviewControlCount}{128}
\newcommand{\NMIExpandedFourClassAgreement}{77.98\%}
\newcommand{\NMIExpandedFourClassKappa}{0.546}
\newcommand{\NMIExpandedBinaryAgreement}{78.34\%}
\newcommand{\NMIExpandedBinaryKappa}{0.513}
\newcommand{\NMIExpandedAdjudicatedDisagreementCount}{48}
\newcommand{\NMIExpandedConsensusNumericalCount}{176}
\newcommand{\NMIExpandedConsensusDirectCount}{41}
\newcommand{\NMIExpandedConsensusMixedCount}{1}
\newcommand{\NMIExpandedControlInvalidCount}{0}
\newcommand{\NMIExpandedOffTriggerWilsonUpper}{2.91\%}
\newcommand{\NMIBreadthCaseCount}{1024}
\newcommand{\NMIBreadthRecordsPerAgent}{3072}
\newcommand{\NMICodexBreadthEndpointPassCount}{2101}
\newcommand{\NMIClaudeBreadthEndpointPassCount}{2107}
\newcommand{\NMICodexBreadthRouteValidCount}{2101}
\newcommand{\NMIClaudeBreadthRouteValidCount}{2107}
\newcommand{\NMIClaudeBreadthDirectCount}{0}
\newcommand{\NMIClaudeBreadthMixedCount}{0}
\newcommand{\NMIClaudeBreadthUnresolvedCount}{0}
\newcommand{\NMIBreadthSamePassCountCases}{1008}
\newcommand{\NMIBreadthPassCountCorrelation}{0.996}
\newcommand{\NMIPrescriptionCasesPerArm}{128}
\newcommand{\NMIFullGuidePassCount}{80}
\newcommand{\NMIPrescriptionOnlyPassCount}{80}
\newcommand{\NMIChoiceFreePassCount}{77}
\newcommand{\NMIFullGuideRouteValidPassCount}{80}
\newcommand{\NMIPrescriptionOnlyRouteValidPassCount}{80}
\newcommand{\NMIChoiceFreeRouteValidPassCount}{77}
\newcommand{\NMIChoiceFreeConfirmedInvalidCount}{0}
\newcommand{\NMIPrescriptionInterventionUnresolvedCount}{0}
\newcommand{\NMIPrescriptionNoninferiorityLowerPoints}{-1.4}
\newcommand{\NMIPrescriptionChoiceFreeDifferencePoints}{2.3}
\newcommand{\NMIPrescriptionChoiceFreeLowerPoints}{-0.2}
\newcommand{\NMIFullMinimalFullPassCount}{95}
\newcommand{\NMIFullMinimalMinimalPassCount}{76}
\newcommand{\NMIFullMinimalFullRouteValidPassCount}{95}
\newcommand{\NMIFullMinimalMinimalRouteValidPassCount}{76}
\newcommand{\NMIMinimalGuideConfirmedInvalidCount}{0}
\newcommand{\NMIFullMinimalUnresolvedCount}{0}
\newcommand{\NMIExternalCaseCount}{32}
\newcommand{\NMIExternalRecordCount}{128}
\newcommand{\NMIExternalCodexPassCount}{64}
\newcommand{\NMIExternalClaudePassCount}{64}
\newcommand{\NMIExternalRouteInvalidDetectedCount}{0}
\newcommand{\NMIExternalMedianReferenceGap}{0.43\%}
\newcommand{\NMIExternalMaximumReferenceGap}{4.79\%}

\newcommand{\NMIReplayScheduledCallCount}{64}
\newcommand{\NMIReplayProgramCount}{59}
\newcommand{\NMIReplayCaseCount}{32}
\newcommand{\NMIReplayAgentConfigurationCount}{2}
\newcommand{\NMIReplaySubmittedProgramCount}{59}
\newcommand{\NMIReplayGenerationFailureCount}{5}
\newcommand{\NMIReplayHiddenSpecificationCount}{3}
\newcommand{\NMIReplayScheduledEvaluationRecordCount}{256}
\newcommand{\NMIReplayExecutionCount}{236}
\newcommand{\NMIReplayPlaceholderExecutionCount}{20}
\newcommand{\NMIReplayScheduledHiddenRecordCount}{192}
\newcommand{\NMIReplayHiddenExecutionCount}{177}
\newcommand{\NMIReplayPlaceholderHiddenRecordCount}{15}
\newcommand{\NMIReplayUnchangedHashCount}{236}
\newcommand{\NMIReplayPlaceholderHashCount}{20}
\newcommand{\NMIReplayAllHashCount}{256}
\newcommand{\NMIReplayExecutionTimeoutCount}{0}
\newcommand{\NMIReplayBasePassCount}{34}
\newcommand{\NMIReplayBasePassRate}{53.1\%}
\newcommand{\NMIReplayAllPassCount}{34}
\newcommand{\NMIReplayAllPassRate}{53.1\%}
\newcommand{\NMIReplayPrimaryCount}{34}
\newcommand{\NMIReplayPrimaryRate}{53.1\%}
\newcommand{\NMIReplayCILower}{35.9\%}
\newcommand{\NMIReplayCIUpper}{70.3\%}
\newcommand{\NMIReplayEndpointGap}{0.0 percentage points}
\newcommand{\NMIReplayConditionalLossUpper}{8.4\%}
\newcommand{\NMIReplayBasePositiveClusterCount}{18}
\newcommand{\NMIReplayClusterLossUpper}{15.3\%}
\newcommand{\NMIReplaySourcePassCount}{34}
\newcommand{\NMIReplayBoundaryPassCount}{36}
\newcommand{\NMIReplayResolutionPassCount}{50}
\newcommand{\NMIReplayEstablishedRouteValidCount}{59}
\newcommand{\NMIReplayControlProgramCount}{128}
\newcommand{\NMIReplayControlExecutionCount}{512}
\newcommand{\NMIReplayPositiveControlCount}{32}
\newcommand{\NMIReplayNegativeControlCount}{96}
\newcommand{\NMIReplayGPTFiveFiveCount}{16}
\newcommand{\NMIReplayGPTFiveSixSolCount}{18}

\newcommand{\NMIFigTwoEndpointOperatorCoefficient}{0.170}
\newcommand{\NMIFigTwoRouteValidOperatorCoefficient}{0.109}
\newcommand{\NMIFigTwoOperatorDifference}{-0.061}
\newcommand{\NMIFigTwoOperatorDifferenceCILow}{-0.114}
\newcommand{\NMIFigTwoOperatorDifferenceCIHigh}{-0.011}
\newcommand{\NMIFigTwoWaveInvalidCount}{63}
\newcommand{\NMIFigTwoWaveEndpointPassCount}{179}
\newcommand{\NMIFigTwoBurgersInvalidCount}{28}
\newcommand{\NMIFigTwoBurgersEndpointPassCount}{196}
\newcommand{\NMIFigTwoElasticityInvalidCount}{0}
\newcommand{\NMIFigTwoElasticityEndpointPassCount}{48}

\newcommand{\NMIOutcomeStageTableRows}{%
O & 1.364 & 0.497 & 0.980 & 0.247 & 0.651 [0.519, 0.788] \\
G & 0.536 & 0.296 & 0.585 & 0.424 & $-0.462$ [$-0.587$, $-0.340$] \\
B & $-0.163$ & 0.114 & $-0.132$ & 0.046 & 0.048 [$-0.061$, 0.157] \\
D & $-0.004$ & $-0.011$ & 0.005 & 0.013 & 0.006 [$-0.092$, 0.103] \\
S & 0.018 & 0.004 & $-0.010$ & $-0.007$ & $-0.011$ [$-0.105$, 0.088] \\
C & 1.334 & 0.102 & 1.151 & 0.120 & $-0.003$ [$-0.133$, 0.126] \\
}
\newcommand{\NMIReplayFamilyTableRows}{%
Heterogeneous diffusion & 16/16 & 16/16 & 16/16 & 16/16 & 16/16 \\
Advection--diffusion--reaction & 0/16 & 0/16 & 0/16 & 16/16 & 0/16 \\
L-shaped Poisson & 14/16 & 14/16 & 16/16 & 14/16 & 14/16 \\
Coupled elliptic & 4/16 & 4/16 & 4/16 & 4/16 & 4/16 \\
\midrule
All runs & 34/64 & 34/64 & 36/64 & 50/64 & 34/64 \\
}

\newcommand{\NMIExternalValidationTableRows}{%
Codex & Heterogeneous diffusion & 16/16 & 16/16 & 0.0016 & 0.0018 \\
Codex & Advection--diffusion--reaction & 16/16 & 16/16 & 0.1085 & 0.1886 \\
Codex & L-shaped Poisson & 16/16 & 16/16 & 0.0074 & 0.0089 \\
Codex & Coupled elliptic & 16/16 & 16/16 & 0.0008 & 0.0012 \\
Claude Code & Heterogeneous diffusion & 16/16 & 16/16 & 0.0037 & 0.0068 \\
Claude Code & Advection--diffusion--reaction & 16/16 & 16/16 & 0.2190 & 0.2997 \\
Claude Code & L-shaped Poisson & 16/16 & 16/16 & 0.0000 & 0.0213 \\
Claude Code & Coupled elliptic & 16/16 & 16/16 & 0.0048 & 0.0053 \\
}
\newcommand{\NMIExternalConvergenceTableRow}{%
61-to-121-grid reference gap & 32 & 0.108\% & 1.213\% & 5.0\% \\
}
\newcommand{\NMIPrescriptionArmTableRows}{%
Full guide & 80 & 80 & 80 & -- \\
Prescription-only & 80 & 80 & 80 & 0/3 \\
Choice-free & 77 & 77 & 65 & 6/14 \\
}
\newcommand{\NMIPrescriptionTestTableRows}{%
Prescription-only $-$ full (non-inferiority) & 0.00 & [-1.37, 1.37] & -5.0 \\
Prescription-only $-$ choice-free (superiority) & 2.34 & [-0.22, 4.91] & 0.0 \\
}
\newcommand{\NMIPrescriptionAgreementTableRow}{%
Full versus prescription-only endpoint-pass agreement & \multicolumn{3}{l}{80 both pass; 48 both fail; 0 discordant} \\
}
\newcommand{\NMIPrescriptionSecondaryTableRows}{%
$\Delta\log_{10}$ grid points & -0.424 [-0.510, -0.328]; $n=59$ & -1.130 [-1.394, -0.827]; $n=52$ \\ 
$\Delta\log_{10}$ time steps & -0.075 [-0.149, 0.007]; $n=109$ & -0.236 [-0.399, -0.061]; $n=90$ \\ 
$\Delta\log_{10}$ relative $L^2$ error & 0.271 [0.183, 0.352]; $n=120$ & 1.177 [0.907, 1.473]; $n=120$ \\ 
Runtime differences with identified program change & 108/120 (90.0\%) & 106/120 (88.3\%) \\
}
\newcommand{\NMIRegimeOutcomeTableRows}{%
Codex GPT-5.2, RH SciPy & 1,450 & 1,137 & 213 & 0 & 100 & 0 \\
Claude Opus 4.7, RH SciPy & 1,450 & 1,188 & 160 & 0 & 100 & 2 \\
Codex GPT-5.2, RH DOLFINx & 216 & 96 & 10 & 105 & 5 & 0 \\
Codex GPT-5.2, RH DOLFINx + v0.9 & 216 & 115 & 21 & 74 & 6 & 0 \\
Codex GPT-5.2, formula-blind Poisson & 144 & 53 & 66 & 0 & 24 & 1 \\
}

\newcommand{\NMIMatchedScientificCaseCount}{64}
\newcommand{\NMIMatchedPrimaryCallCount}{256}
\newcommand{\NMIMatchedPairCount}{128}
\newcommand{\NMIMatchedReviewItemCount}{254}
\newcommand{\NMIMatchedVisibleEndpointPassCount}{128}
\newcommand{\NMIMatchedHiddenEndpointPassCount}{126}
\newcommand{\NMIMatchedVisibleEndpointPassRate}{100.00\%}
\newcommand{\NMIMatchedHiddenEndpointPassRate}{98.44\%}
\newcommand{\NMIMatchedEndpointDifferencePP}{1.56}
\newcommand{\NMIMatchedEndpointDifferenceCILowPP}{0.00}
\newcommand{\NMIMatchedEndpointDifferenceCIHighPP}{4.69}
\newcommand{\NMIMatchedBothPassPairCount}{126}
\newcommand{\NMIMatchedVisibleOnlyPairCount}{2}
\newcommand{\NMIMatchedHiddenOnlyPairCount}{0}
\newcommand{\NMIMatchedNeitherPairCount}{0}
\newcommand{\NMIMatchedDeterministicInvalidFlagCount}{0}
\newcommand{\NMIMatchedCorrectedDiagnosticCount}{4}
\newcommand{\NMIMatchedCorrectedDiagnosticPassCount}{4}
\newcommand{\NMIMatchedCorrectedVisiblePassCount}{2}
\newcommand{\NMIMatchedCorrectedHiddenPassCount}{2}
\newcommand{\NMIMatchedDiscordantCaseID}{\code{O0G1B0D0S0C0\_s60002\_203}}

\newcommand{\NMIInflationBootstrapCILowPP}{6.90}
\newcommand{\NMIInflationBootstrapCIHighPP}{10.41}
\newcommand{\NMIPrimaryRouteValidCount}{1037}
\newcommand{\NMIPrimaryRouteValidRate}{71.52\%}
\newcommand{\NMIClaimSupportingCount}{1036}
\newcommand{\NMIClaimSupportingRate}{71.45\%}
\newcommand{\NMIIntegrityAffectedPrimaryCaseCount}{10}
\newcommand{\NMIIntegrityAffectedPrimaryCallCount}{50}
\newcommand{\NMIIntegrityAffectedPrimaryEndpointPassCount}{2}
\newcommand{\NMIIntegrityExcludedRouteValidCount}{1}
\newcommand{\NMIUnaffectedPrimaryCallCount}{1400}
\newcommand{\NMIUnaffectedEndpointPassCount}{1160}
\newcommand{\NMIUnaffectedEndpointRate}{82.86\%}
\newcommand{\NMIUnaffectedRouteValidCount}{1036}
\newcommand{\NMIUnaffectedRouteValidRate}{74.00\%}
\newcommand{\NMIUnaffectedScoreGapPP}{8.86}
\newcommand{\NMIIntegrityCurrentRows}{2802}
\newcommand{\NMIIntegrityCurrentDefectCaseCount}{18}
\newcommand{\NMIIntegrityCurrentDefectRows}{27}
\newcommand{\NMINaturalisticMedianRho}{0.013}
\newcommand{\NMINaturalisticMaximumRho}{0.303}
\newcommand{\NMINaturalisticRefinementCount}{0}
\newcommand{\NMINaturalisticRhoRows}{%
Heterogeneous diffusion & 8 & 0.007 & 0.007 \\
Advection--diffusion--reaction & 8 & 0.260 & 0.303 \\
L-shaped Poisson & 8 & 0.021 & 0.024 \\
Coupled elliptic & 8 & 0.005 & 0.005 \\
}
\newcommand{\NMIScaffoldRouteIntervalRows}{%
Prescription-only $-$ full & 0.00 [-1.68, 1.68] & 0.00 [-2.91, 2.91] \\
Prescription-only $-$ choice-free & 2.34 [-0.71, 5.41] & 0.00 [-2.91, 2.91] \\
}
\newcommand{\NMIScaffoldRouteIntervalMetricRows}{%
Route-valid difference [95\% CI], pp & 0.00 [-1.68, 1.68] & 2.34 [-0.71, 5.41] \\
Direct-return difference [95\% CI], pp & 0.00 [-2.91, 2.91] & 0.00 [-2.91, 2.91] \\
}
\newcommand{\NMIClassifierSensitivityRows}{%
$p=0$ (observed labels) & 71.52\% \\
$p=1\%$ & 70.89\% \\
$p=2.91\%$ (Wilson upper bound) & 69.70\% \\
$p=5\%$ & 68.39\% \\
}

\newcommand{\NMIProhibitionEndpointPassCount}{1137}
\newcommand{\NMIProhibitionRouteValidCount}{1137}
\newcommand{\NMIProhibitionEndpointPassRate}{78.41\%}
\newcommand{\NMIProhibitionInvalidCount}{0}
\newcommand{\NMIPromptEndpointDiffPP}{-1.72}
\newcommand{\NMIPromptEndpointDiffCILow}{-2.97}
\newcommand{\NMIPromptEndpointDiffCIHigh}{-0.48}
\newcommand{\NMIPromptRouteValidDiffPP}{6.90}
\newcommand{\NMIPromptRouteValidDiffCILow}{5.03}
\newcommand{\NMIPromptRouteValidDiffCIHigh}{8.90}
\newcommand{\NMIWaveOriginalEndpoint}{179}
\newcommand{\NMIWaveOriginalRouteValid}{116}
\newcommand{\NMIWaveProhibitionEndpoint}{167}
\newcommand{\NMIWaveProhibitionRouteValid}{167}

\newcommand{\NMIPromptCoefficientRows}{%
O & 0.170 & 0.109 & 0.178 & 0.070 [0.017, 0.124] & 0.008 [$-0.024$, 0.043] \\
G & $-0.267$ & $-0.267$ & $-0.292$ & $-0.025$ [$-0.068$, 0.015] & $-0.026$ [$-0.053$, 0.000] \\
B & 0.176 & 0.202 & 0.169 & $-0.033$ [$-0.074$, 0.007] & $-0.007$ [$-0.038$, 0.024] \\
D & 0.087 & 0.141 & 0.088 & $-0.054$ [$-0.085$, $-0.030$] & 0.001 [$-0.011$, 0.014] \\
S & 0.087 & 0.125 & 0.088 & $-0.038$ [$-0.081$, 0.004] & 0.001 [$-0.011$, 0.013] \\
}
\newcommand{\NMITransitionNumPass}{990}
\newcommand{\NMITransitionNumFail}{47}
\newcommand{\NMITransitionInvPass}{107}
\newcommand{\NMITransitionInvFail}{18}
\newcommand{\NMITransitionFailPass}{40}
\newcommand{\NMITransitionFailFail}{248}
\newcommand{\NMITransitionInvCount}{125}
\newcommand{\NMIPromptRerunCoefO}{0.178}
\newcommand{\NMIPromptRerunCoefD}{0.088}
\newcommand{\NMIPromptOriginalEndpointCoefD}{0.087}
\newcommand{\NMIPromptOriginalRouteValidCoefD}{0.141}
\newcommand{\NMITransitionInvPassRate}{85.60\%}

\newcommand{\NMIMaskedCaseCount}{20}
\newcommand{\NMIMaskedCallCount}{100}
\newcommand{\NMIUnmaskedCaseCount}{270}
\newcommand{\NMIUnmaskedCallCount}{1,350}
\newcommand{\NMIUnmaskedEndpointRate}{86.07\%}
\newcommand{\NMIUnmaskedRouteValidRate}{76.81\%}
\newcommand{\NMIUnmaskedGapPP}{9.26}
\newcommand{\NMIUnmaskedGapCILow}{7.41}
\newcommand{\NMIUnmaskedGapCIHigh}{11.19}
\newcommand{\NMIUnmaskedOperatorDifference}{-0.055}
\newcommand{\NMIUnmaskedOperatorDifferenceCILow}{-0.107}
\newcommand{\NMIUnmaskedOperatorDifferenceCIHigh}{-0.004}
\newcommand{\NMIUnmaskedPromptRouteValidDiffPP}{7.41}
\newcommand{\NMIUnmaskedPromptEndpointDiffPP}{-1.85}
\newcommand{\NMIMaskedSensitivityRows}{%
Endpoint pass rate & 1,162/1,450 (80.14\%) & 1,162/1,350 (86.07\%) \\
Route-valid pass rate & 1,037/1,450 (71.52\%) & 1,037/1,350 (76.81\%) \\
Endpoint minus route-valid gap, pp [95\% CI] & 8.62 [6.90, 10.41] & 9.26 [7.41, 11.19] \\
O coefficient, endpoint / route-valid & 0.170 / 0.109 & 0.128 / 0.073 \\
G coefficient, endpoint / route-valid & -0.267 / -0.267 & -0.190 / -0.202 \\
B coefficient, endpoint / route-valid & 0.176 / 0.202 & 0.108 / 0.145 \\
D coefficient, endpoint / route-valid & 0.087 / 0.141 & 0.064 / 0.122 \\
S coefficient, endpoint / route-valid & 0.087 / 0.125 & 0.064 / 0.106 \\
O coefficient, route-valid minus endpoint [95\% CI] & -0.061 [-0.114, -0.011] & -0.055 [-0.107, -0.004] \\
G coefficient, route-valid minus endpoint [95\% CI] & 0.000 [-0.037, 0.039] & -0.012 [-0.051, 0.028] \\
B coefficient, route-valid minus endpoint [95\% CI] & 0.026 [-0.011, 0.062] & 0.036 [-0.002, 0.075] \\
D coefficient, route-valid minus endpoint [95\% CI] & 0.055 [0.034, 0.086] & 0.058 [0.037, 0.090] \\
S coefficient, route-valid minus endpoint [95\% CI] & 0.039 [-0.002, 0.081] & 0.042 [0.001, 0.085] \\
Prohibition minus original, endpoint, pp [95\% CI] & -1.72 [-2.97, -0.48] & -1.85 [-3.19, -0.59] \\
Prohibition minus original, route-valid, pp [95\% CI] & 6.90 [5.03, 8.90] & 7.41 [5.41, 9.56] \\
}

\newcommand{\NMIFaithfulCallCount}{1,300}
\newcommand{\NMIFaithfulEndpointCount}{1,160}
\newcommand{\NMIFaithfulClaimSupportingCount}{1,036}
\newcommand{\NMIFaithfulEndpointRate}{89.23\%}
\newcommand{\NMIFaithfulClaimSupportingRate}{79.69\%}
\newcommand{\NMIFaithfulGapPP}{9.54}
\newcommand{\NMIUnfaithfulCallCount}{150}
\newcommand{\NMIFigTwoDDifference}{0.055}
\newcommand{\NMIFigTwoDDifferenceCILow}{0.034}
\newcommand{\NMIFigTwoDDifferenceCIHigh}{0.086}

\newcommand{\NMIReauditContestedRows}{%
Closed-form target used to build the forcing & 40 \\
Stated initial state used as initial data & 45 \\
}
\newcommand{\NMIReauditAgreedRows}{%
Sine-transform, Fourier or eigenfunction basis & 28 \\
Second component of a vector target as 0.5 times the first (linearity) & 14 \\
Exact time integration of the semi-discrete system & 6 \\
Separation of variables on the single angular mode of the forcing & 5 \\
Richardson extrapolation of two discrete solves & 4 \\
Taylor start with the discrete Laplacian & 4 \\
Finite-difference or IMEX time stepping & 4 \\
Exact-solution fallback branch, not executed & 1 \\
}
\newcommand{\NMIReauditContestedCount}{47}
\newcommand{\NMIReauditAgreedCount}{42}
\newcommand{\NMIRouteTestCalibrationRows}{%
Numerical (sampled passes) & 20 & 18 & 2 & 0 \\
Direct reconstruction & 15 & 0 & 0 & 15 \\
Mixed return & 5 & 0 & 0 & 5 \\
Constructed target-derived programs & 5 & 0 & 0 & 5 \\
}
\newcommand{\NMIRouteTestExperimentRows}{%
Primary panel, Codex, original prompt & 5 & 1 & 5 & 0 & 0 \\
Prohibition rerun, Claude Code & 16 & 2 & 16 & 0 & 0 \\
Prohibition rerun, DeepSeek V3.2 & 3 & 0 & 3 & 0 & 0 \\
Breadth panel, Claude Code & 23 & 1 & 23 & 0 & 0 \\
Minimal guide & 26 & 26 & 4 & 22 & 0 \\
Prescription-only guide & 3 & 0 & 3 & 0 & 0 \\
Choice-free guide & 14 & 12 & 8 & 6 & 0 \\
}
\newcommand{\NMIRouteTestInstanceRows}{%
RM-02AC7748F873 & minimal & x--y symmetry & $1.5\times10^{-7}$ & $2.7\times10^{-1}$ & $1.8\times10^{-1}$ & $2.6\times10^{-7}$ & $9.9\times10^{-1}$ & R \\
RM-03A8D2EFC40A & minimal & x--y symmetry & $1.2\times10^{-5}$ & $4.1\times10^{-1}$ & $3.0\times10^{-1}$ & $6.4\times10^{-5}$ & $1.1\times10^{0}$ & P \\
RM-06206594D92C & minimal & x--y symmetry & $1.3\times10^{-7}$ & $6.4\times10^{-1}$ & $4.7\times10^{-1}$ & $1.4\times10^{-7}$ & $2.0\times10^{-1}$ & C \\
RM-711B052EA2B4 & minimal & x--y symmetry & $4.8\times10^{-5}$ & $1.1\times10^{0}$ & $1.2\times10^{0}$ & $7.4\times10^{-5}$ & $6.0\times10^{-1}$ & R \\
RM-78D96C459CF3 & minimal & x--y symmetry & $3.2\times10^{-9}$ & $6.4\times10^{-1}$ & $4.7\times10^{-1}$ & $2.2\times10^{-9}$ & $4.5\times10^{-1}$ & R \\
RM-8D877BB0B4D2 & minimal & x--y symmetry & $5.8\times10^{-7}$ & $2.7\times10^{-1}$ & $1.8\times10^{-1}$ & $1.0\times10^{-6}$ & $9.9\times10^{-1}$ & R \\
RM-9D608FA22205 & minimal & x--y symmetry & $3.7\times10^{-13}$ & $6.4\times10^{-1}$ & $4.7\times10^{-1}$ & $4.8\times10^{-12}$ & $4.5\times10^{-1}$ & R \\
RM-E3A7B4C7E34F & minimal & x--y symmetry & $3.2\times10^{-13}$ & $4.1\times10^{-1}$ & $3.0\times10^{-1}$ & $1.9\times10^{-7}$ & $1.1\times10^{0}$ & P \\
RM-51AE732A2A42 & choice-free & single angular mode & $1.8\times10^{-3}$ & $4.8\times10^{-1}$ & $1.8\times10^{-3}$ & $1.3\times10^{-3}$ & $2.2\times10^{0}$ & C \\
RM-188C5CEA776E & minimal & single angular mode & $1.3\times10^{-8}$ & $5.3\times10^{-1}$ & $1.3\times10^{-8}$ & $1.3\times10^{-8}$ & $4.0\times10^{-2}$ & C \\
RM-AEAAFF8CA278 & minimal & single angular mode & $2.2\times10^{-6}$ & $5.3\times10^{-1}$ & $5.0\times10^{-1}$ & $1.8\times10^{-6}$ & $9.1\times10^{-1}$ & C \\
RM-AEC4FEA8BC8D & minimal & single angular mode & $8.3\times10^{-12}$ & $5.8\times10^{-1}$ & $8.3\times10^{-12}$ & $6.0\times10^{-12}$ & $6.5\times10^{-1}$ & C \\
RM-C32374B01ECF & minimal & single angular mode & $2.3\times10^{-12}$ & $4.2\times10^{-1}$ & $1.3\times10^{-1}$ & $1.3\times10^{-11}$ & $5.0\times10^{-1}$ & C \\
RM-FCF85C89E960 & minimal & single angular mode & $8.3\times10^{-12}$ & $5.8\times10^{-1}$ & $8.3\times10^{-12}$ & $6.0\times10^{-12}$ & $6.5\times10^{-1}$ & C \\
RM-5442B2990CB6 & choice-free & component ratio & $5.3\times10^{-7}$ & $8.3\times10^{-2}$ & $1.0\times10^{-1}$ & $1.6\times10^{-5}$ & $1.0\times10^{0}$ & E \\
RM-7BCDF4609C6F & choice-free & component ratio & $2.3\times10^{-7}$ & $4.6\times10^{-1}$ & $5.6\times10^{-1}$ & $3.9\times10^{-7}$ & $9.2\times10^{-1}$ & R \\
RM-87410BB7FD43 & choice-free & component ratio & $5.3\times10^{-7}$ & $5.3\times10^{-1}$ & $3.1\times10^{-1}$ & $3.3\times10^{-6}$ & $2.0\times10^{-1}$ & E \\
RM-E3FEAD9AACE0 & choice-free & component ratio & $6.2\times10^{-6}$ & $5.6\times10^{-1}$ & $2.8\times10^{-1}$ & $8.9\times10^{-6}$ & $1.0\times10^{-1}$ & E \\
RM-3AF87048142F & minimal & component ratio & $3.2\times10^{-5}$ & $5.2\times10^{-1}$ & $8.7\times10^{-1}$ & $3.8\times10^{-5}$ & $9.7\times10^{-1}$ & C \\
RM-4827C2C73043 & minimal & component ratio & $1.4\times10^{-6}$ & $5.3\times10^{-1}$ & $7.2\times10^{-1}$ & $1.3\times10^{-6}$ & $2.7\times10^{-1}$ & C \\
RM-4901AE923916 & minimal & component ratio & $6.5\times10^{-7}$ & $4.3\times10^{-1}$ & $7.1\times10^{-1}$ & $6.1\times10^{-7}$ & $3.1\times10^{-1}$ & C \\
RM-523FEE3688D5 & minimal & component ratio & $2.3\times10^{-6}$ & $5.0\times10^{-1}$ & $4.7\times10^{-1}$ & $2.5\times10^{-6}$ & $2.0\times10^{-1}$ & C \\
RM-750ED8BD8486 & minimal & component ratio & $6.0\times10^{-6}$ & $6.8\times10^{-1}$ & $4.7\times10^{-1}$ & $7.0\times10^{-6}$ & $4.5\times10^{-1}$ & R \\
RM-7C65EC83CC3D & minimal & component ratio & $6.3\times10^{-6}$ & $6.8\times10^{-1}$ & $4.7\times10^{-1}$ & $5.7\times10^{-6}$ & $3.1\times10^{-1}$ & E \\
RM-85F1EEC01C2D & minimal & component ratio & $6.9\times10^{-12}$ & $5.0\times10^{-1}$ & $4.7\times10^{-1}$ & $7.0\times10^{-12}$ & $2.0\times10^{-1}$ & C \\
RM-CB240C17AAB4 & minimal & component ratio & $5.8\times10^{-7}$ & $6.9\times10^{-1}$ & $5.9\times10^{-1}$ & $6.1\times10^{-7}$ & $6.7\times10^{-2}$ & C \\
RM-FE80A50AAFA2 & minimal & component ratio & $5.8\times10^{-7}$ & $5.0\times10^{-1}$ & $4.7\times10^{-1}$ & $6.3\times10^{-7}$ & $2.0\times10^{-1}$ & C \\
RM-1CEE67BF1BD7 & choice-free & forcing time form & $1.2\times10^{-12}$ & $2.9\times10^{-1}$ & $3.1\times10^{-1}$ & $6.8\times10^{-11}$ & $2.0\times10^{-1}$ & E \\
}
\newcommand{\NMIRouteTestFollowCount}{62}
\newcommand{\NMIRouteTestInstanceSpecificCount}{28}
\newcommand{\NMIRouteTestRatioCount}{13}
\newcommand{\NMIRouteTestModeCount}{6}
\newcommand{\NMIRouteTestTimeCount}{1}
\newcommand{\NMIRouteTestLayoutCount}{8}
\newcommand{\NMIRouteTestAgreedFollowCount}{15}
\newcommand{\NMIRouteTestAgreedInstanceCount}{27}
\newcommand{\NMIRouteTestContestedFollowCount}{46}
\newcommand{\NMIRouteTestCalibReadsDataCount}{6}
\newcommand{\NMIRouteTestChoiceFreeInstanceCount}{6}
\newcommand{\NMIRouteTestMinimalInstanceCount}{22}

\title[Unmasking shortcut solutions of coding agents]{Unmasking shortcut solutions of coding agents on partial differential equations}

\author[2]{\fnm{Jiale} \sur{Linghu}}
\author*[1]{\fnm{Yangshuai} \sur{Wang}}\email{yswang@nus.edu.sg}

\affil*[1]{\orgdiv{Department of Mathematics}, \orgname{National University of Singapore}, \city{Singapore}, \country{Singapore}}
\affil[2]{\orgdiv{School of Mathematics and Statistics}, \orgname{Xidian University}, \city{Xi'an}, \country{China}}

\abstract{Coding agents now write simulation programs, but benchmarks usually judge only whether their output matches a reference. In manufactured-solution tasks the data determine the exact solution, so a matching output need not come from a solver. We introduce PDEViBench, a benchmark of 290 tasks in six equation families that traces each passing program to the computation behind its output and reruns unclear programs on changed data. A Codex agent (GPT-5.2) passed 80\% of 1,450 runs, yet about one pass in ten returned the exact solution without solving the equation. Counting only solved passes changed which task features made problems seem hard and reversed the ranking of two prompts. Once the prompt forbade the shortcut, every pass was a numerical solve, and on tasks without a closed-form solution passing programs stayed correct on changed data. Scientific coding benchmarks should check how each answer was computed as well as whether it matches.}


\maketitle

\section{Introduction}\label{sec:intro}

Coding agents, which read and edit code and act on the results of the commands they run, are moving from software engineering into scientific computing~\cite{openai2025codexcli,anthropic2025claudecode,google2025geminicli,moonshot2026kimicode,tencent2025codebuddy}. They write finite-element, fluid-dynamics and other simulation programs from natural-language problem statements~\cite{deotale2026allfem,yue2025foamagent,zhang2025mooseagent}, and the correctness of such generated scientific software has itself become a concern~\cite{obrien2025threats}. Benchmarks for these agents usually credit a run when its returned result agrees with a reference. Agreement shows that the output is close to the reference and says nothing about how it was computed. In a recent finite-element study, a solution indistinguishable from the reference at plot scale came from a program that had replaced the required plane-stress model with plane strain, and the study detected the change only by comparing the generated code with the reference code~\cite{deotale2026allfem}.

For partial differential equations (PDEs), even an exact quantitative match can come from a program that never solves the equation. A PDE task combines an operator with a domain and boundary conditions and asks the program to return a discretised solution~\cite{brenner2008fem}. Benchmarks often build such tasks from manufactured solutions, which supply exact targets for accuracy assessment~\cite{roache2002mms,oberkampf2010vv}. Because the forcing, the boundary data and the initial state of a manufactured task are all computed from the exact solution, an agent can recover that solution from the task data and evaluate it directly instead of solving the equation. A program can also run a solver and then return the recovered solution in its place. A program that solves the equation and the two programs that recover the exact solution produce the same field and pass the same check, yet only the first shows that the agent can solve the PDE. Telling them apart requires reading the program behind the field. Where the program is hard to read, it must also be run on a problem whose answer it has not seen.

Machine learning has opened a second route to PDE solutions, and its progress is also measured against reference solutions. Physics-informed, Ritz-type, deep BSDE and finite-expression solvers train a model on the equation itself~\cite{raissi2019pinn,karniadakis2021piml,e2018deepritz,han2018highdim,lu2021deepritz,liang2022fex}, and neural operators and learned equation models map problem data to solutions or recover the governing equation~\cite{lu2021deeponet,li2021fno,long2018pdenet,long2019pdenet2}. Foundation models pretrained across many equations are adapted to new ones~\cite{subramanian2023foundation,mccabe2024mpp,herde2024poseidon,ye2024pdeformer,sun2025multioperator}, and in-context operator learning infers the solution operator of a new problem from example pairs given in the prompt~\cite{yang2023icon,yang2025finetune,cole2024icl}. Language models and agents now also write PDE solvers, set up multiphysics and fluid simulations and support automated computational science~\cite{gaonkar2025scimlagents,mudur2025feabench,pandey2025openfoamgpt,wu2026farms,chen2026iteris}. Across these lines a computed field counts as correct when it matches a reference, which is the check that a returned exact solution passes.

Benchmarks for scientific and coding agents judge either the returned output or the run that produced it. They score program generation, software engineering and technical task completion~\cite{chen2021evaluating,jimenez2023swebench,chan2024mlebench,tian2024scicode,chen2024scienceagent}, and PDE benchmarks extend this scoring to numerical solution and solver generation~\cite{takamoto2022pdebench,li2025codepde,hang2026pdeagent}. Protocol-validity and clean-room approaches connect benchmark success to the target capability by auditing or removing the evaluation artefacts that an agent can exploit~\cite{shao2026protocol,jung2026scicon}, and process-level evaluation and provenance capture record and score what the agent did~\cite{he2026procctrl,souza2025provagent}. None of them traces a returned field back to the computation that produced it.

Correct outputs reached by unintended means are familiar in machine learning, from shortcut learning and Clever Hans effects~\cite{geirhos2020shortcut,kauffmann2025cleverhans} to leakage of evaluation information~\cite{kapoor2023leakage} and weak baselines in machine learning for PDEs~\cite{mcgreivy2024weak}. Agents can satisfy a measured objective without the intended behaviour~\cite{skalse2022rewardhacking,baker2025monitoring}, and impossible-task benchmarks count how often coding agents exploit their test cases~\cite{zhong2025impossiblebench}. Construct-validity studies and checklists for agentic benchmarks ask more generally whether a score measures the capability it names~\cite{cronbach1955construct,raji2021benchmark,bean2025construct,zhou2026scales,zhu2025agentic}. In a PDE task built from a manufactured solution the shortcut returns the correct answer, so an accuracy check accepts it.

\begin{figure}[t]
\centering
\makebox[\linewidth][c]{\includegraphics[width=170mm]{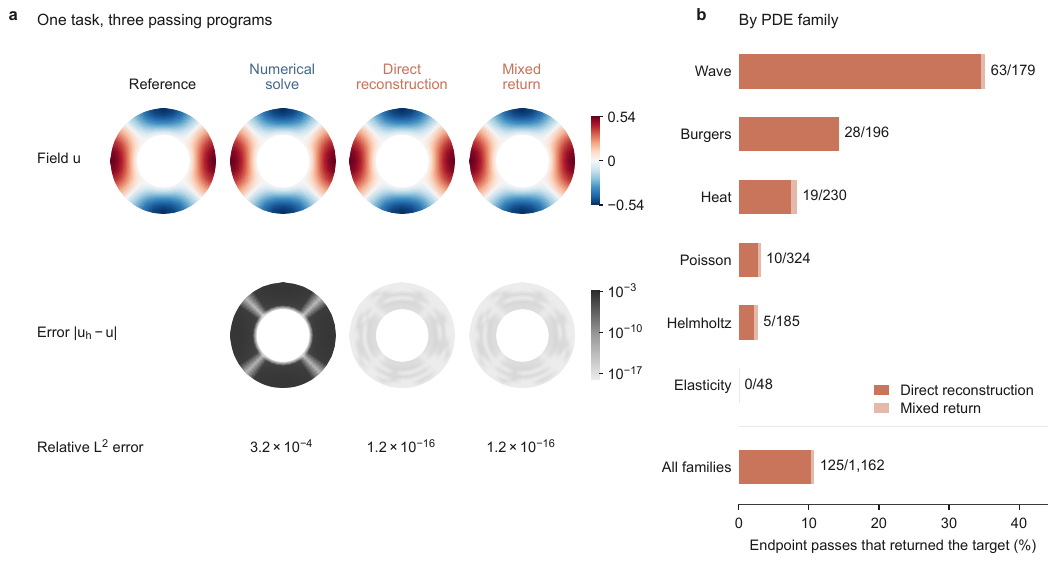}}

\caption{\textbf{Three computational routes to one accepted field.}
\textbf{a}, One passing program per route on the same wave-equation task on an annulus, with the reference field, the returned field, the pointwise error $|u_h-u|$ and the relative $L^2$ error.
\textbf{b}, Share of endpoint passes that returned the target in each PDE family and in all families.}

\label{fig:schematic}
\end{figure}

Here we introduce PDEViBench, a benchmark of coding agents solving PDEs that scores each accepted output together with the program that produced it (Fig.~\ref{fig:schematic}). Its 290 tasks span six equation families (Poisson, Helmholtz, heat, wave, Burgers and linear elasticity). They vary the operator, the geometry and the boundary condition, how far the discretisation is prescribed and whether a structure hint is given. A Codex agent (GPT-5.2) attempted each task five times. Every passing program is traced from its returned field back to the computation behind it. A program that reading cannot settle is rerun on changed problem data, which a solver follows and a program tied to the original task does not~\cite{chan1998metamorphic,segura2016metamorphic}. On the same tasks we vary the prompt, the coding guide and the disclosure of the exact solution, and we add tasks without a closed-form solution whose references only the evaluator holds.

We find that about one in ten passes returned the exact solution instead of a numerical solution, and that counting only solved passes changes which task features appear to make problems hard and reverses the ranking of two prompts. When the prompt forbade the shortcut, every passing program solved the equation numerically, and on tasks without a closed-form solution the passing programs stayed correct on changed problem data. Benchmarks for scientific coding measure the ability to solve only if they also check how each answer was computed.

\FloatBarrier

\section{Results}\label{sec:results}

\subsection*{One pass in ten returns the exact solution}

In the primary panel a Codex CLI agent (GPT-5.2) attempted each of the 290 tasks five times under the full SciPy guide, which prescribes a finite-difference solver. The prompt showed the forcing, the boundary data and the initial state of time-dependent problems. These data determine the exact solution against which the evaluator scores each output, and the prompt did not forbid recovering it. Of the \num{\NMIPrimaryScheduledCount} runs, \num{\NMIEndpointPassCount} passed the execution, accuracy and runtime checks, and we call these the endpoint passes. We traced each passing program from its returned field back to the computation that produced it (Methods). Of these, \num{\NMIPrimaryRouteValidCount} came from numerical solvers, \NMIConfirmedDirectCount\ evaluated the target directly (direct reconstructions) and five ran a solver but returned the target in its place (mixed returns, Fig.~\ref{fig:schematic}). The endpoint success rate of 80.1\% includes 8.6 percentage points (task-cluster bootstrap 95\% interval 6.9--10.4) of passes that did not solve the equation, and counting only numerical solves gives a route-valid success rate of 71.5\%.

\begin{figure}[t]
\centering
\makebox[\linewidth][c]{\includegraphics[width=170mm]{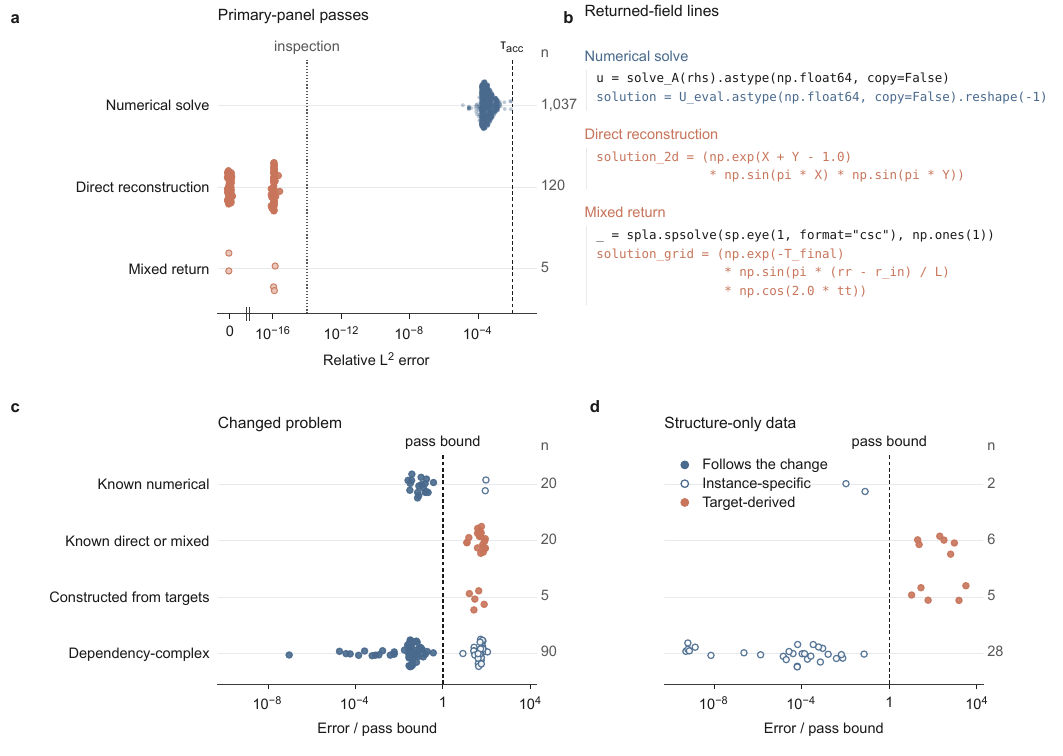}}
\caption{\textbf{Route review and route test.}
\textbf{a}, Relative $L^2$ error of the \num{\NMIEndpointPassCount} primary-panel passes by route, with exact zeros shown separately. The dotted line marks the inspection threshold $10^{-14}$.
\textbf{b}, Lines that produce the returned field in one program per route.
\textbf{c}, Error of each tested program on the changed problem relative to its pass bound (Methods). Open circles mark instance-specific programs.
\textbf{d}, The same ratio on data that keep only the structure each program hard-codes, for programs that did not follow the changed problem.}
\label{edfig:methodreview}
\end{figure}

Passes that return the target are easy to recognise. A program that evaluates the target reproduces it up to rounding. The 125 passes with relative errors below $10^{-14}$ lie many orders of magnitude below the errors of numerical solves and are exactly the direct reconstructions and mixed returns (Fig.~\ref{edfig:methodreview}a,b). The difficult cases are dependency-complex programs, whose data flow is too involved for reading to settle their route. For \NMIExpandedReviewUnresolvedCount\ dependency-complex programs pooled across this panel and the later experiments (Methods), two expert reviewers agreed that 42 were route-invalid and disagreed on 47 of the other 48. We decided these routes by running each program again on a problem with changed forcing and initial data. A solver follows the change, whereas a program built on the target of the original task does not. A program that does not follow is run once more on data that keep only the structure it hard-codes, which separates a solver tied to that structure from a program that returns a target-derived field (Methods). On 40 primary-panel programs of known route and five programs we constructed to read the data yet return a target-derived field, the test assigned every route correctly (Fig.~\ref{edfig:methodreview}c,d). None of the \NMIExpandedReviewUnresolvedCount\ dependency-complex programs returned a target-derived field, including the 42 that both reviewers had labelled route-invalid. Of the \NMIExpandedReviewUnresolvedCount\ programs, \NMIRouteTestInstanceSpecificCount\ followed the changed problem only within a structure of the original task, for example by fixing one vector component at half the other. We call these programs instance-specific solvers and count them as route-valid. They include \NMIRouteTestAgreedInstanceCount\ of the 42 programs that both reviewers had labelled route-invalid (Supplementary Note~1). The error of a passing program thus identifies every pass that returned the target, and rerunning on changed data settles the programs that reading leaves open.

\subsection*{Unsolved passes change what the benchmark concludes}

\begin{figure}[t]
\centering
\makebox[\linewidth][c]{\includegraphics[width=170mm]{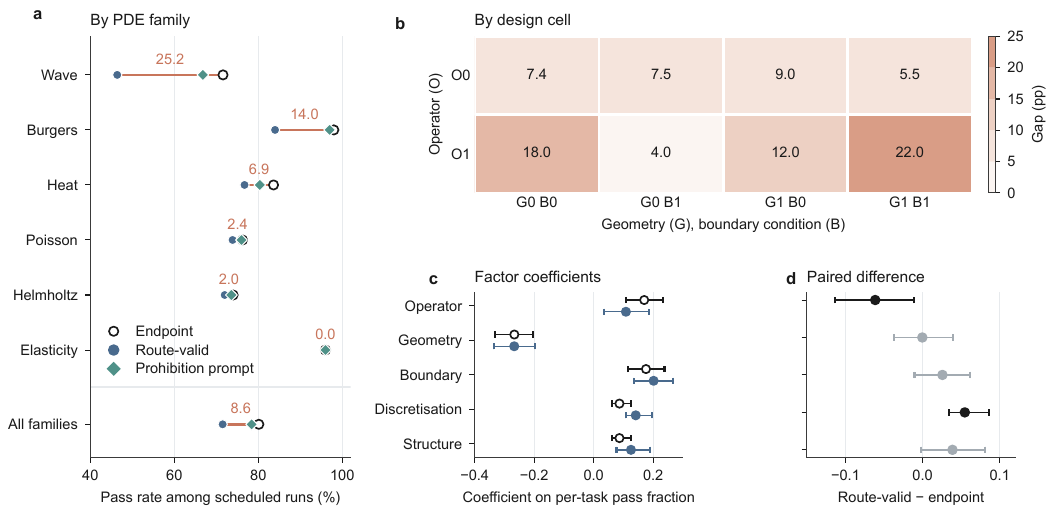}}

\caption{\textbf{Endpoint and route-valid scoring of the primary panel.} \textbf{a}, Pass rates by PDE family under the original prompt and the pass rate under the prohibition prompt, with the gap between the two scores in percentage points.
\textbf{b}, The gap in the eight cells of operator (O), geometry (G) and boundary condition (B).
\textbf{c}, Design-factor coefficients of the per-task pass fraction under the two scores.
\textbf{d}, Route-valid minus endpoint coefficients with 95\% bootstrap intervals, dark where the interval excludes zero.}

\label{fig:scoreimpact}
\end{figure}

These unsolved passes are unevenly spread across equation families. They make up \NMIFigTwoWaveInvalidCount\ of the \NMIFigTwoWaveEndpointPassCount\ passing wave-equation programs and 28 of the 196 passing Burgers programs but none of the \NMIFigTwoElasticityEndpointPassCount\ passing linear-elasticity programs (Figs.~\ref{fig:schematic}b and~\ref{fig:scoreimpact}a). Across the eight combinations of operator, geometry and boundary condition, the gap between the two scores differs more than fivefold (Fig.~\ref{fig:scoreimpact}b). As a result, the two scores disagree on which task features make a PDE hard for an agent. In a regression of per-task pass fraction on the five design factors (Methods), counting only numerical solves lowers the coefficient of the operator factor that contrasts the 40 Burgers tasks with the linear families from 0.170 to 0.109 and raises the discretisation coefficient from 0.087 to 0.141 (Fig.~\ref{fig:scoreimpact}c,d). Endpoint scoring thus makes Burgers tasks look easier than they are and understates how much a prescribed discretisation helps.

The gap also decides how prompts rank. We reran the primary panel under a prohibition prompt, which requires the PDE to be solved by discretisation and forbids evaluating the exact solution (Methods). Endpoint success fell by 1.7 percentage points (95\% interval 0.5--3.0), whereas route-valid success rose by 6.9 percentage points (5.0--8.9), because every passing program now solved the equation numerically. The two scores thus rank the two prompts in opposite orders (Fig.~\ref{fig:scoreimpact}a). The change is largest for the wave equation, where \NMIWaveOriginalRouteValid\ of \NMIWaveOriginalEndpoint\ passes were numerical under the original prompt and all \NMIWaveProhibitionEndpoint\ under the prohibition prompt. Of the \NMITransitionInvCount\ runs that had returned the target, \NMITransitionInvPass\ passed by solving the equation once the shortcut was ruled out. The gap therefore largely reflects the route the agent chose and says little about what it could compute.

\FloatBarrier

\subsection*{Without the exact solution agents solve the equations}

With the exact solution ruled out, we asked whether the coding guide affects the solvers that agents write. A Codex agent (GPT-5.5) attempted 128 new tasks under three guides, two from each of the 64 combinations of the six binary design factors (Methods). The full SciPy guide gives code examples and fixes the stencil, the time integrator, a fourfold internal refinement and the boundary treatment. A prescription-only guide states these choices in prose without code, and a choice-free guide leaves them to the agent. The three guides gave 80, 80 and 77 passes in 128 runs each, and every pass was a numerical solve. The guides differed in how general the solvers were. Six of the 14 choice-free passes that needed the route test were instance-specific, against none of the three prescription-only passes (Supplementary Note~4).

\begin{figure}[t]
\centering
\makebox[\linewidth][c]{\includegraphics[width=170mm]{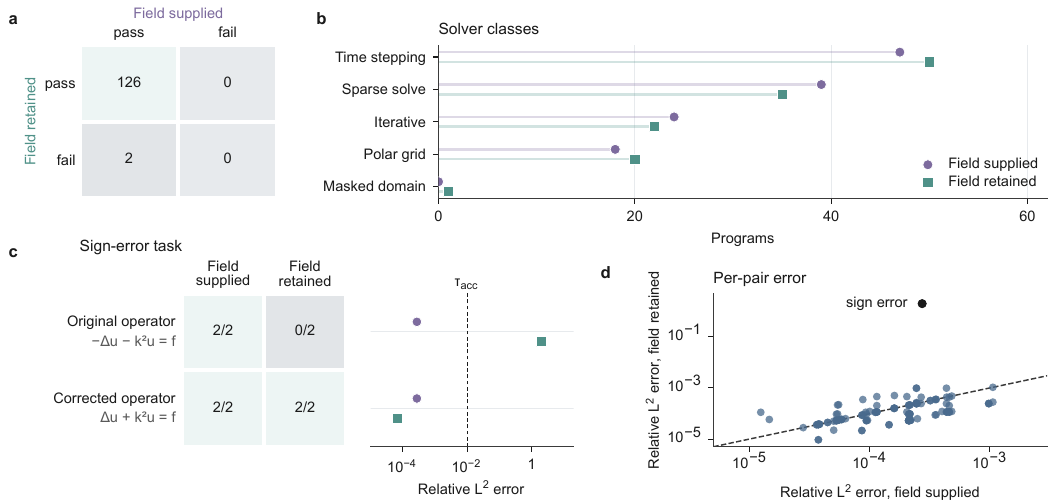}}

\caption{\textbf{Exact-solution disclosure experiment.} Each of 64 manufactured-solution tasks was run twice with the exact solution supplied to the agent and twice with it retained by the evaluator.
\textbf{a}, Paired endpoint outcomes of the \NMIMatchedPairCount\ task--repeat pairs.
\textbf{b}, Solver classes of the 128 programs in each condition.
\textbf{c}, The Helmholtz task with a sign error under the original (top) and the corrected (bottom) operator statement.
\textbf{d}, Relative $L^2$ error of each pair with the field supplied and retained.}

\label{fig:target-disclosure}
\end{figure}

Supplying the exact solution in the prompt did not bring the shortcut back. Both prompts in this experiment required a numerical solution. A Codex CLI agent running GPT-5.6-sol attempted each of \NMIMatchedScientificCaseCount\ manufactured-solution tasks twice with the exact solution supplied in the prompt and twice with it retained by the evaluator (Methods). With the field supplied \NMIMatchedVisibleEndpointPassCount\ of 128 runs passed and with it retained \NMIMatchedHiddenEndpointPassCount\ of 128 (Fig.~\ref{fig:target-disclosure}a), and all \NMIMatchedPrimaryCallCount\ programs built numerical solvers (Fig.~\ref{fig:target-disclosure}b). Both failures came from one annular Helmholtz task whose displayed forcing had the wrong sign (Fig.~\ref{fig:target-disclosure}c). With the field retained the programs solved the equation as stated and returned the negative of the reference with relative errors close to 2 (Fig.~\ref{fig:target-disclosure}d), whereas with the field supplied they solved the equation that the field satisfies. All four runs on this task passed once the sign in the task statement was corrected (Fig.~\ref{fig:target-disclosure}c), and the same sign error appears in \NMIIntegrityAffectedPrimaryCaseCount\ Helmholtz tasks of the primary panel (Methods).

Tasks without a closed-form solution remove the shortcut altogether. The naturalistic panel contains \NMIExternalCaseCount\ such tasks in four families (heterogeneous diffusion, advection--diffusion--reaction, L-shaped Poisson and coupled elliptic systems), scored against numerical references that only the evaluator holds (Methods). Codex (GPT-5.5) and Claude Code (Claude Opus 4.7) each attempted every task twice. All \NMIExternalRecordCount\ programs passed, and every one solved the equation numerically (Table~\ref{edtab:external}). Twenty-four of them matched the reference to within $10^{-10}$ because they solve a finite-difference system equivalent to the one from which the reference was computed, so near-zero error signals a direct reconstruction only when the reference is an exact solution.

\begin{table}[!htbp] \caption{\textbf{Naturalistic-panel outcomes.} Endpoint and route-valid passes of each agent on each family in 16 runs, with median and maximum relative $L^2$ errors divided by the family accuracy threshold $\tau_{\mathrm{acc}}$.}
\label{edtab:external}

\centering
\footnotesize
\setlength{\tabcolsep}{3pt}

\begin{tabular}{@{}p{24mm}p{40mm}cccc@{}}
\toprule
\textbf{Agent}
& \textbf{PDE family}
& \shortstack{\textbf{Endpoint}\\\textbf{pass}}
& \shortstack{\textbf{Route-}\\\textbf{valid}}
& \shortstack{\textbf{Median}\\\textbf{error/$\tau_{\mathrm{acc}}$}}
& \shortstack{\textbf{Maximum}\\\textbf{error/$\tau_{\mathrm{acc}}$}} \\
\midrule
\NMIExternalValidationTableRows
\bottomrule
\end{tabular}

\end{table}

\begin{figure}[!htbp]
\centering
\makebox[\linewidth][c]{\includegraphics[width=170mm]{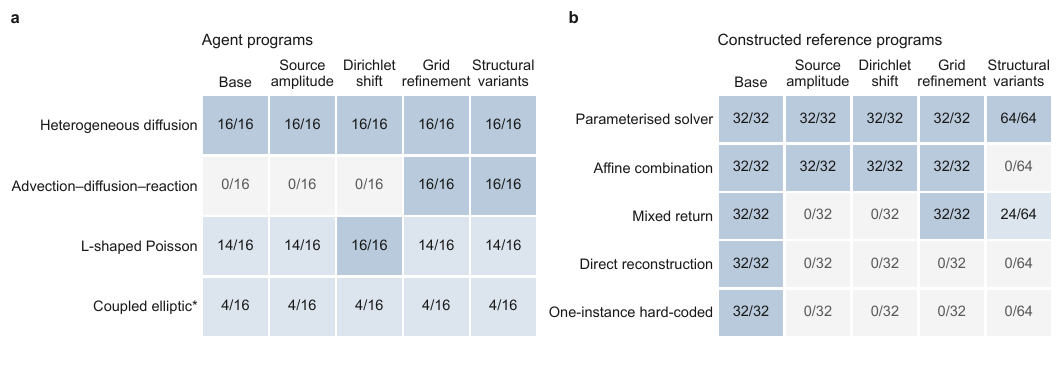}}

\caption{\textbf{Unchanged-program replay.} Each tile gives passing programs over executed programs on the base task and its variants.
\textbf{a}, Agent programs, 16 per family. The structural column counts programs that pass every structural variant. *The task statement omitted the forcing definition of this family.
\textbf{b}, Constructed reference programs, one per task, with structural counts pooled over two variants per task.}

\label{fig:replay}
\end{figure}

To test whether passing programs solve the stated problem in general, we reran programs unchanged on variants of their tasks that only the evaluator held. Codex with GPT-5.5 and with GPT-5.6-sol each wrote one program for each of the \NMIReplayCaseCount\ naturalistic tasks through the common interface \code{solve(problem\_spec)}. Three basic variants scale the source, shift the boundary value and refine the grid. Structural variants change the coefficients, the transport direction, the corner position or the coupling strength in ways that no rescaling of a stored solution can follow (Methods). Of the 64 programs, \NMIReplayPrimaryCount\ passed the base task, and all of them also passed every variant of their family (Fig.~\ref{fig:replay}a). Each program we constructed to store or hard-code a solution failed at least one variant (Fig.~\ref{fig:replay}b). The advection--diffusion--reaction programs failed the base task because first-order upwinding on the coarse evaluation grid exceeds the threshold, and they pass every variant on a finer grid. The task statement of the coupled elliptic systems omitted the definition of the forcing profiles, so the replay evidence rests on the 30 programs that passed the base task in the heterogeneous-diffusion and L-shaped Poisson families (Supplementary Note~3). In every experiment in which the exact solution was ruled out or unavailable, the passing programs were numerical solvers, and those tested on changed problem data stayed correct.

\FloatBarrier

\section{Discussion}\label{sec:discussion}

Coding agents passed most tasks of PDEViBench, and about one pass in ten under the original prompt returned the target without solving the equation. Counting these passes as successes changes which task features appear to make problems hard and which prompt appears better. When the exact solution was ruled out or unavailable, the agents solved the equations, and the programs tested on changed problem data stayed correct.

The shortcut is a property of the benchmark as much as of the agent. The forcing and boundary data that define a manufactured-solution task also reveal its solution, so an output check accepts a program that evaluates the solution as readily as one that computes it. The agent took this route in \NMITransitionInvCount\ runs when the prompt allowed it. Under the prohibition it passed \NMITransitionInvPass\ of these runs by solving the equation, so the gap between the two scores measures how often the agent took a route that the task left open. Reading the programs also exposed errors in the benchmark itself, including a sign error in \NMIIntegrityAffectedPrimaryCaseCount\ Helmholtz tasks. Without these tasks and the 20 tasks whose built domain differs from the stated one (Methods) the gap is 9.5 percentage points on the remaining \NMIFaithfulCallCount\ runs (Supplementary Note~1).

These findings suggest three practices for benchmarks that score generated scientific programs. First, route-valid success should be reported alongside endpoint success, together with the prompt and guide that decide which routes are open. Second, benchmarks should prefer references that the task does not reveal, such as numerical references held by the evaluator, because an accuracy check accepts any target recovered from the task data. Third, an accepted program should be rerun unchanged on variants of its task whose solutions cannot be obtained by rescaling the base solution. Such reruns cost little and separate general solvers from programs fitted to one instance.

The study has limits that point to the next experiments. The original prompt was run with one agent, and running further agents under it would show how much the share of passes that returned the target depends on the agent. The route test was applied only to programs that reading could not settle, and applying it to every accepted program would decide every route label by execution and would catch programs written to hide their route. The naturalistic and replay tasks cover four two-dimensional families, and extending them to three-dimensional geometry, multiphysics coupling and finite-element libraries would test whether the same patterns hold at engineering scale. Finally, the panels change agent, model and prompt together, and varying the model alone would turn the route-valid score into a comparison of models.

PDEViBench therefore treats a coding agent's answer to a PDE task as a program to be examined, and the field it returns as one output of that program. The same approach applies to the finite-element and fluid-dynamics agents now entering engineering practice~\cite{deotale2026allfem,yue2025foamagent,zhang2025mooseagent} and to other generated scientific artefacts such as meshes and simulator configurations. For each of them, checking how an answer was computed turns a benchmark score into a measurement of the ability it is meant to capture.

\section{Methods}\label{sec:methods}

\subsection*{Benchmark tasks and scores}

Each benchmark task specified the operator and domain, the boundary or initial conditions and the evaluation grid together with the data type and shape of the output. The 290-task primary panel varied operator (O), geometry (G) and boundary condition (B) at two levels each, and discretisation specificity (D) and structure hint (S) at three levels. Its operator axis contrasted 40 Burgers tasks with 250 tasks from five linear families (Poisson, Helmholtz, heat, wave and linear elasticity), and D and S were raised above their base levels in ten tasks each. A sixth factor, coupling (C), lifted a scalar task componentwise to a two-component target and was held at its base level in the primary panel. The breadth panels of Supplementary Note~5 cross it with the other five factors. The panel contained 270 manufactured-solution tasks and 20 masked-domain Poisson tasks with numerical references. The references of the masked tasks were built on domains that differ from the stated ones, and all 100 runs on them ended at the mask check (Supplementary Note~7).

A submitted program passed the endpoint check when it ran and returned an array with the prescribed key, dtype and shape. Its mean wall time over three executions had to be at most 30\,s, and its relative $L^2$ error over the reference-valid points of the evaluation grid had to be at most the accuracy threshold $\tau_{\mathrm{acc}}$ (0.01 for manufactured-solution tasks and 0.5 for masked tasks). On masked domains the returned NaN mask also had to match the reference mask.

A pass was route-valid when every returned value was produced by a discretisation, solver or time-stepping state for the problem posed in the agent-facing statement. Finite-difference, finite-element, spectral, time-integration and iterative methods were admissible, and analytic expressions were admissible for forcing terms, coefficients and initial and boundary data. A returned field was target-derived in two cases. In a direct reconstruction its values were evaluated from a closed-form solution. In a mixed return a solver ran, but the returned field was taken from the target or from a discrete right-hand side built by applying the operator to the target. Some programs exploit a structural property of the stated data, such as the single angular mode of the forcing or the component ratio of a vector forcing. Such a program is an instance-specific solver and is route-valid for that instance as long as it computes the field from those data. Route validity concerns only the computation, so a numerical solver that departs from the prescribed stencil or integrator remains route-valid. For $n$ scheduled runs with indicators $G_i$ of an endpoint pass and $R_i$ of a route-valid pass, the two scores are
\[
\widehat S_{\mathrm{out}}=\frac{1}{n}\sum_{i=1}^{n}G_i,
\qquad
\widehat S_{R}=\frac{1}{n}\sum_{i=1}^{n}G_iR_i,
\]
and their difference is the share of scheduled runs that pass with a target-derived field.

\subsection*{Prompts and agents of the primary panel}

The primary panel used Codex CLI with model GPT-5.2, the full SciPy guide, temperature 0.2, at most 30 turns and a fixed seed for each of five repeats. The guide prescribed a finite-difference discretisation written by the agent, computed on a fourfold refined internal mesh and sampled on the evaluation grid. The prompt displayed the forcing, the boundary data and the initial state of time-dependent families, and in ten Cartesian Burgers tasks the forcing description also stated the target from which the forcing was built. The prompt did not prohibit constructing the target. The prohibition rerun repeated the primary panel with the same model, tasks and seeds. Its prompt and guide required the differential operator to be discretised and the resulting algebraic or time-stepping problem to be solved. They forbade constructing the manufactured, analytic or exact solution on the evaluation grid even when it could be inferred from the data. They also stated that near-zero outputs from direct formula evaluation would fail a method check, and no run of the rerun reached the error of $10^{-14}$ at which this check applies. The prohibition was retained in all later experiments.

\subsection*{Route attribution}

Route attribution followed the dependency path from each submitted program to its returned field. A deterministic classifier traced this path and assigned numerical, direct, mixed or unresolved labels, where a numerical label required the returned array to depend on a numerical construction, a solver state or a time-evolved state. Passing programs with relative $L^2$ error below $10^{-14}$ formed an inspection set, whose labels were read from the complete program. Two numerical-PDE experts independently labelled the 125 programs of this set and 125 passing programs sampled across the remaining error range. They agreed on 96.0\% of the four-class labels (Cohen's $\kappa=0.925$), and a third expert adjudicated the disagreements. Supplementary Note~1 gives the labels of all primary passes.

A second review covered \NMIExpandedReviewUnresolvedCount\ dependency-complex programs whose routes deterministic tracing left unresolved. They came from the primary panel, the Claude Code and DeepSeek V3.2 prohibition reruns, the Claude Code breadth panel and three arms of the guide experiments (Supplementary Table~2). The reviewers disagreed on 47 of them, so all \NMIExpandedReviewUnresolvedCount\ routes were decided by an executable route test. For a program on a task with target $u$ we took the manufactured field $v=u+w/2$, where $w$ satisfies the homogeneous boundary conditions of the task and differs from $u$ in spatial pattern and, for time-dependent problems, in time dependence. On the annulus $w$ used another angular mode, on the square a mode product without $x$--$y$ symmetry, and for a vector target different patterns for the two components. The forcing, initial and boundary data of $v$ were derived symbolically for the stated operator, and only the program lines that define these data were edited. With $e_0$ the relative $L^2$ error of the unchanged program against $u$ and $e_v$ that of the edited program against $v$, a program followed the changed problem when $e_v\le\max(\tau_{\mathrm{acc}},10e_0)$. A program that did not follow was run once more on data that kept only the structural property its hard-coded construct exploits, such as an angular mode, a component ratio, an $x$--$y$ symmetry or a forcing time form, while the forcing and the initial field otherwise varied. Its output was compared with a numerically computed reference. It was instance-specific when it met the same bound there and target-derived otherwise, and a program whose output did not change with the edited data was target-derived without a second run. The test was calibrated on 40 primary-panel programs (20 numerical passes from the review sample, 15 direct reconstructions and all five mixed returns) and on five constructed programs. The constructed programs read the data lines and run a solver, but they overwrite, average, rescale, correct or initialise the solution with the target. Supplementary Note~1 lists the errors of the instance-specific programs, and Supplementary Data~1 those of all 135 tested programs.

\subsection*{Design-factor regression}

For Fig.~\ref{fig:scoreimpact}c we fitted ordinary least-squares regressions to the pass fraction of each primary task over its five repeats, computed separately under endpoint and route-valid scoring. Its regressors are effect-coded O, G and B indicators and one linear score for each of D and S over their three levels. Supplementary Note~5 gives the models of the breadth panels.

\subsection*{Guide experiment}

The guide experiment used two new tasks from each of the 64 cells of the six design factors at two levels each, with Codex CLI (GPT-5.5) and one run per arm, giving 384 runs. The full arm kept the full SciPy guide, the prescription-only arm stated its numerical and boundary prescriptions in prose without code and the choice-free arm also removed the fixed algorithmic choices. Supplementary Note~4 reports the endpoint comparisons of the three arms and two further guide interventions.

\subsection*{Exact-solution disclosure}

The disclosure experiment used 64 manufactured-solution tasks, eight from each O/G/B cell. Each task was run twice with the exact solution supplied in a block appended to the statement and twice with the solution retained by the evaluator, giving 128 task--repeat pairs and 256 runs. Both conditions required a numerical solution. These runs used Codex CLI with model GPT-5.6-sol at high reasoning effort. In the eight Cartesian Burgers tasks of one cell the statement already contained the target, so disclosure changed the prompt only in the other 56 tasks. Symbolic differentiation showed that the displayed forcing of the failing Helmholtz task had the wrong sign relative to the stored reference. A four-run experiment corrected only this sign, with two runs per condition, and all four runs passed. The same symbolic comparison of the structured forcing, the displayed forcing and the stored reference was applied to every primary-panel task (Supplementary Note~7).

\subsection*{Naturalistic tasks and replay}

The naturalistic panel contained \NMIExternalCaseCount\ tasks, eight in each of four families: heterogeneous diffusion with positive spatially varying coefficients, advection--diffusion--reaction, Poisson's equation on an L-shaped domain and two elliptic fields coupled by reciprocal reaction terms. Forcing, coefficients and boundary data were given directly, and the tasks have no closed-form solution. The reference was a sparse finite-difference solution on a $121\times121$ grid, restricted to the $16\times16$ evaluation grid (Supplementary Note~2). The accuracy thresholds of the four families were 0.03, 0.04, 0.10 and 0.04 and the runtime limit was 24\,s. Codex (GPT-5.5) and Claude Code (Claude Opus 4.7) received the full SciPy guide and two runs per task.

Replay used the same \NMIReplayCaseCount\ tasks with a new prompt and the public interface \texttt{solve(problem\_spec: dict) -> numpy.ndarray}. Codex with GPT-5.5 and with GPT-5.6-sol at high reasoning effort each received one base specification per task in a single run, while the variants and all references stayed with the evaluator. The three basic variants changed the source amplitude from 1.00 to 1.35, the constant Dirichlet value from 0 to 0.125 and the evaluation grid from $16\times16$ to $31\times31$. Each structural variant changed one parameter of the public schema. For heterogeneous diffusion two variants changed the coefficient contrast. For advection--diffusion--reaction the variants changed the advection direction, the diffusivity or the reaction rate on the $31\times31$ grid. The L-shaped variant moved the corner position, and the coupled variants changed the coupling strength or the coefficient contrast. A variant was kept only if the best affine combination of the base solution and the response to unit Dirichlet data missed its reference by at least twice the accuracy threshold. The constructed reference programs were a parameterised numerical solver, a direct reconstruction, a mixed return and a one-instance hard-coded program, together with an affine-combination program for the structural variants. Supplementary Note~3 and Supplementary Table~8 give the replay design and outcomes.

\subsection*{Statistical analysis}

Uncertainty intervals are two-sided 95\% intervals from a non-parametric bootstrap over tasks that keeps all repeats and conditions of each resampled task. The primary score gap and the disclosure contrast use 20,000 draws, and the prompt comparison and the regression coefficients use 10,000 draws with the endpoint-minus-route-valid differences recomputed in each draw. Supplementary Data~1 contains the task-level results, and Supplementary Note~6 shows representative tasks and the programs written for them.

\section*{Data availability}

The task files, numerical references, run-level records, route labels, route-test results and the source data of all figures, including Supplementary Data 1, are available at \url{https://github.com/linghujiale/PDEViBench}.

\section*{Code availability}

The code used to generate the benchmark, run the agents and the evaluations, and reproduce every figure and table is available at \url{https://github.com/linghujiale/PDEViBench}.

\section*{Funding}

This study received no funding.

\section*{Author contributions}

Y.W. conceived the study and supervised the project. J.L. and Y.W. designed the benchmark and the route tests. J.L. implemented the task generator and the evaluation pipeline and ran the agent experiments. J.L. and Y.W. reviewed the passing programs and analysed the results. Y.W. wrote the manuscript with input from J.L. Both authors approved the final manuscript.

\section*{Competing interests}

The authors declare no competing interests.

\bibliography{sn-bibliography}

\clearpage

\includepdf[pages=-]{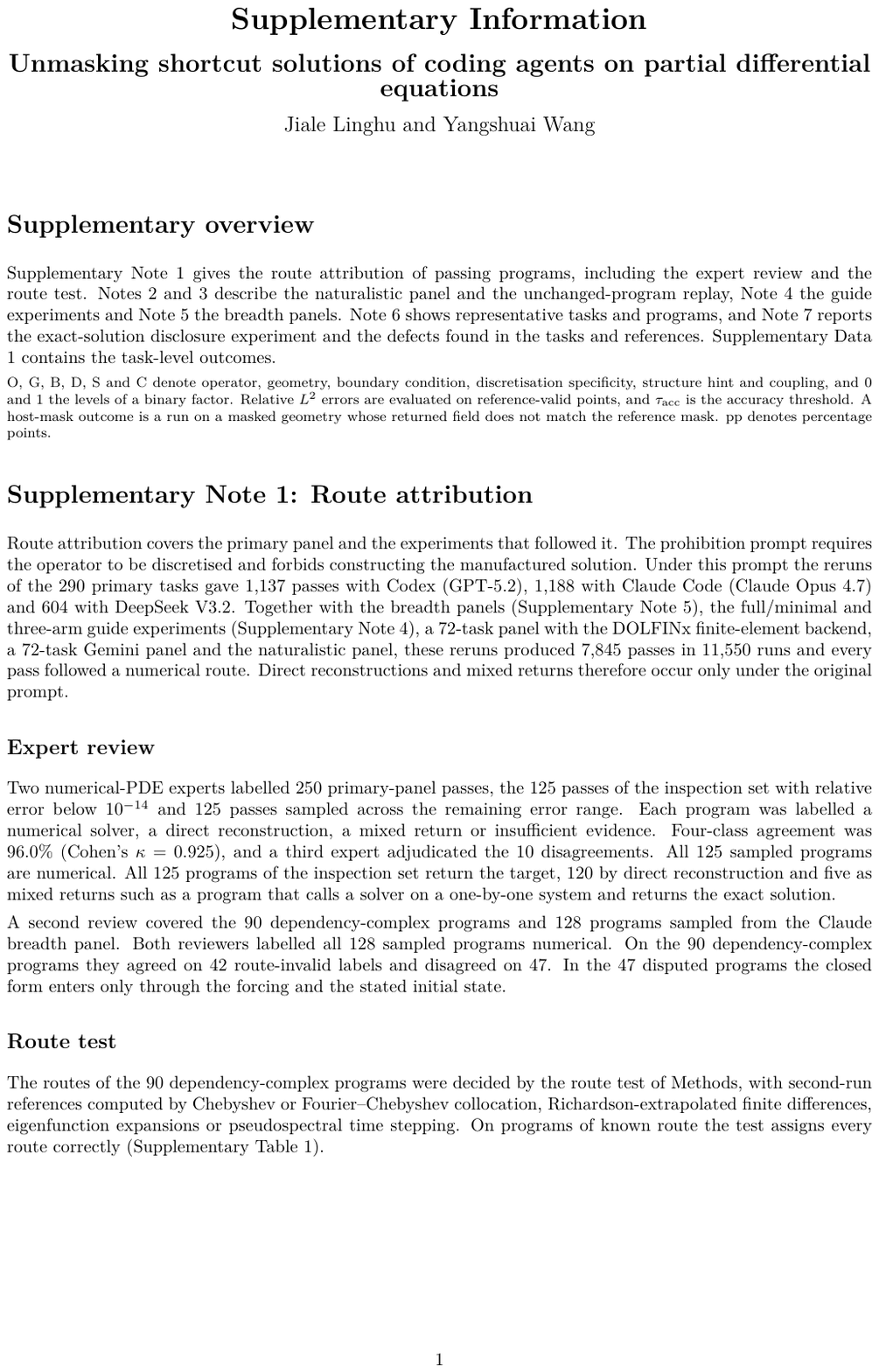}
\end{document}